\documentclass[12pt]{article}

\usepackage{amsmath,amsthm,amssymb,amsopn,amsfonts,pdfpages,url}
\usepackage[linesnumbered,ruled,vlined]{algorithm2e}
\usepackage{graphics,relsize}
\usepackage{graphicx}
\usepackage{multirow}
\usepackage{subfig}
\usepackage{bbm}
\usepackage{authblk}
\usepackage[square,numbers]{natbib}

\newcommand{\twospace}{\renewcommand{\baselinestretch}{1.3}\normalsize}

\newcommand{\p}{\mathbb{P}}
\newcommand{\q}{\mathbbm{1}}
\newcommand{\e}{\mathbb{E}}

\newcommand{\spairs}{{[n] \choose 2}}
\newcommand{\pairs}{\{ i,j \}}

\newcommand{\re}{\varrho_e}
\newcommand{\rh}{\varrho_h}
\newcommand{\rT}{\varrho_T}
\newcommand{\pn}{\varPi_n}
\newcommand{\ppm}{{\mathcal P}_m}
\newcommand{\vc}{\textup{vec}}
\newcommand{\id}{{\mathcal I}}

\newcommand{\st}{\mathfrak{str}}
\newcommand{\dens}{\mathfrak{d}}

\newtheorem{theorem}{Theorem}

\newtheorem{lemma}[theorem]{Lemma}

\twospace
\begin{document}
\title{\vspace{-.8in} Alignment Strength and Correlation for Graphs}
\author[$\dag$]{Donniell E. Fishkind$^\dag$, Lingyao Meng$^\dag$, Ao Sun$^\dag$,\\ Carey E. Priebe}
\author[$\ddag$]{Vince Lyzinski}
\affil[$\dag$]{\small Department of Applied Mathematics and Statistics, Johns Hopkins University\\
Baltimore, MD 21218}
\affil[$\ddag$]{\small Department of Mathematics and Statistics, University of Massachusetts Amherst\\
Amherst, MA 01003}
\maketitle

\begin{abstract}
When two graphs have a correlated Bernoulli distribution, we prove that the alignment strength of their natural bijection
strongly converges to a novel measure of graph correlation $\rT$ that neatly combines intergraph with intragraph distribution
parameters. Within broad families of the random graph parameter settings, we illustrate that exact graph matching runtime
and also matchability are both functions of $\rT$, with thresholding behavior starkly illustrated in matchability.

\footnote{{\tt def@jhu.edu, lmeng2@jhu.edu, asun17@jhu.edu, cep@jhu.edu, vlyzinski@umass.edu}}
\end{abstract}

\newpage

\section{Overview \label{S:overview} } 

Suppose $G$ and $H$ are any two graphs with the same number of vertices.
For any positive integer~$n$, define $[n]:=\{ 1,2,3, \ldots n \}$, and  let
$\spairs$ denote the set of all $2$-element subsets of $[n]$. For simplicity,
suppose that the vertex sets of $G$ and $H$ are both $[n]$.
Let $\pn$ denote the set of bijections from $[n]$ to $[n]$.
For each $\phi \in \pn$, we define the {\it number of disagreements
between $G$ and $H$ under $\phi$} to be
\begin{eqnarray} \label{eqn:astr}
d(G,H,\phi):= \ \ \ \ \ \ \ \ \ \ \ \ \ \ \ \ \ \ \ \ \ \ \ \ \ \ \ \ \ \ \ \ \
\ \ \ \ \ \ \ \ \ \ \ \ \ \ \ \ \ \ \ \ \ \ \ \nonumber \\
\sum_{\pairs \in \spairs}\q \Big ( \ \q \big (i \sim_G j \big )
\ \ne \ \q \big (\phi(i) \sim_H \phi(j) \big ) \ \Big ),
\end{eqnarray}
where $\q (\cdot )$ denotes the indicator function, and
$\sim_G$ denotes adjacency of vertices in $G$.

For each $\phi \in \pn$, we define the {\it alignment strength of $\phi$} as
\begin{eqnarray} \label{eqn:st}
\st (G,H,\phi):=1-\frac{d(G,H,\phi)}{ \frac{1}{n!}\sum_{\phi' \in \pn}d(G,H,\phi')}.
\end{eqnarray}
The denominator in this definition of alignment strength serves as a normalizing factor; in particular, if $\phi$ is an
isomorphism between $G$ and~$H$ then the alignment strength of $\phi$ is $1$,
and if the number of adjacency disagreements for $\phi$ is merely average among the
bijections in $\pn$ then the alignment strength of $\phi$ is $0$.
(Of course, if $G$ and $H$ are
both edgeless or both complete graphs then $\st (G,H,\phi)$ is not defined.)

If $\phi \in \pn$ happens to be a known ``natural alignment" between $G$ and $H$ (for example, if
$G$ and $H$ are social networks with the same members, and $\phi$ maps each member
to themselves; e.g.~an email network and a Twitter network with the same users) then
$\st (G,H,\phi)$ can be viewed as a numerical measure of the structural similarity between $G$ and $H$.
However, if a natural alignment between $G$ and $H$ is not known, then we can use the {\it graph matching problem
solution}, which is defined as $\phi_{GM} \in \arg \min_{\phi' \in \pn}d(G,H,\phi')$; specifically,
$\st (G,H,\phi_{GM})$ can be viewed as a numerical measure of the structural similarity between  $G$ and $H$.

Two practical notes regarding computation:
Although the denominator $ \frac{1}{n!}\sum_{\phi' \in \pn}d(G,H,\phi')$ in the definition of alignment strength (Equation~\ref{eqn:st})
involves an exponentially sized summation, nonetheless it can be computed efficiently using Equation \ref{eqn:b} from Section \ref{S:main}.
Also, although the computation of  the graph matching problem solution~$\phi_{GM}$ is intractable [\cite{thirtyyears}], nonetheless
there are effective, efficient approximate graph matching algorithms that can be used [\cite{VogCon}, \cite{sgm1}], one of which we
discuss and use later in this paper.

A brief outline of this paper is as follows.

In Section \ref{S:randG} we describe a very general random graph setting; $G$ and $H$ are random graphs with a
correlated Bernoulli distribution.
In particular, $G$ and $H$ share the same vertex set, and the identity
bijection $\id \in \pn$ is the natural alignment between $G$ and $H$.
Each pair of vertices is assigned
its own probability of adjacency (``Bernoulli parameter") in $G$ and $H$, and the indicator Bernoulli random
variable for adjacency of the pair in $G$ and the indicator Bernoulli random
variable for adjacency of the pair in $H$ have Pearson correlation coefficient $\re$.
Inherent to this model is the inter-graph (i.e.~between $G$ and $H$) statistical correlation $\re$ and the intra-graph
{\it heterogeneity correlation} parameter $\rh$, which is a function of the Bernoulli coefficients that measures their variation.
Then we define the key parameter $\rT$ as
$1-\rT:=(1-\re)(1-\rh)$; we call $\rT$ the {\it total correlation}.

In Section \ref{S:main} we state and prove our main theoretical result, Theorem \ref{thm:after}, which asserts that
for $G$ and $H$ with a correlated Bernoulli distribution we have that the alignment strength of the identity bijection
$\st(G,H,\id)$ is asymptotically equal to the total correlation parameter $\rT$. This suggests that the total correlation $\rT$ is
a meaningful measure of the structural similarity between the graphs $G$ and $H$ realized from the correlated Bernoulli distribution. Of note is that
the total correlation is nicely and cleanly partitioned by the defining formula  $1-\rT=(1-\re)(1-\rh)$; this illustrates
a symmetry in the affect of (inter-graph parameter) edge correlation $\re$ and the affect
of (intra-graph parameter) heterogeneity correlation $\rh$.

The subsequent sections, Section \ref{S:matchability} and Section \ref{S:runtime}, follow up with empirical illustrations
that total correlation $\rT$ is a meaningful measure. As we vary the edge correlation $\re$ together with the heterogeneity correlation $\rh$
for correlated Bernoulli graphs $G$ and $H$ in broad families of parameter settings, it turns out that the value of $\rT$ dictates
(in Section \ref{S:matchability}) how successful the approximate seeded graph matching algorithm called SGM [\cite{sgm1}, \cite{sgm2}] is in
recovering the identity bijection (which is the natural alignment here) and (in Section \ref{S:runtime}) $\rT$ dictates how much time it
takes to perform seeded graph matching exactly via binary integer linear programming. The {\it seeded graph matching problem} is the graph matching
problem wherein
we seek to compute $\phi_{GM} \in \arg \min_{\phi' \in \pn}d(G,H,\phi')$, except that part of the natural alignment is known;
having these ``seeds" can substantially help recover the rest of the natural alignment correctly. In Section \ref{S:matchability},
we utilize the SGM Algorithm [\cite{sgm1}, \cite{sgm2}] for approximate seeded graph matching on moderately sized graphs, on the order of $1000$ vertices, since,
unfortunately, exact seeded graph matching can only be done on very small, toy-size graphs (a few tens of non-seed vertices).
In Section \ref{S:runtime}, where we are interested in comparing runtime, the approximate seeded graph matching
algorithms are not appropriate to use, since their run times tend to be monolithic (given the number of vertices) and less sensitive
to the parameters of the random graph distribution. So we do exact seeded graph matching, but
only on small~enough~examples.


\section{Random graph setting: correlated Bernoulli graphs \label{S:randG}}

In this section we describe the correlated Bernoulli random graph distribution,
and three important associated parameters/ functions of parameters; namely $\re$, $\rh$, and $\rT$.

For any positive integer $n$, the distribution parameters are any given real number $\re$
(called the {\it edge correlation}) from the interval $[0,1]$, and any given set of real numbers
$\{ p_{i,j} \}_{\pairs \in \spairs}$ (called the {\it Bernoulli parameters}) from the interval
$[0,1]$ such that the Bernoulli parameters are not all equal to $0$ and not all equal to $1$.
Random graphs $G$ and $H$, each on vertex set~$[n]$, will be called
{\it $\re$-correlated random Bernoulli$(\{ p_{i,j} \}_{\pairs \in \spairs})$ graphs}
if, for each $\pairs \in \spairs$, we have that $\q (i \sim_G j)$ is a Bernoulli$(p_{i,j})$ random variable,
and $\q (i \sim_H j)$ is a  Bernoulli$(p_{i,j})$ random variable, and,
if $0<p_{i,j}<1$, then the two random variables $\q (i \sim_G j)$ and $\q (i \sim_H j)$
have Pearson correlation coefficient $\re$; other than these specified
dependencies, the random variables $\{ \q (i \sim_G j)\}_{\pairs \in
\spairs} \bigcup \{ \q (i \sim_H j)\}_{\pairs \in \spairs} $ are collectively independent.

Such $G$, $H$ can be realized from this distribution as follows. For all
$\pairs \in \spairs$ independently, first realize $\q (i \sim_G j)$
from the Bernoulli$(p_{i,j})$ distribution. Then, conditioned on $\q (i \sim_G j)$,
realize $\q (i \sim_H j)$ from distribution Bernoulli$(\re \cdot \q (i \sim_G j) +
(1-\re) \cdot p_{i,j})$. It is easy to verify that $\q (i \sim_H j)$
has a marginal distribution Bernoulli$(p_{i,j})$ and, indeed, the
 random variables $\q (i \sim_G j)$ and $\q (i \sim_H j)$ have Pearson correlation~$\re$ if $0<p_{i,j}<1$.
 Moreover, it easy to verify that, for any two Bernoulli$(p_{i,j})$ random variables such that $0<p_{i,j}<1$,
the Pearson correlation coefficient uniquely determines their joint distribution. Also, it is
easy to verify that $\p [i \sim_G j \ \& \ i \not \sim_H j ]=(1-\re)p_{i,j}(1-p_{i,j})$. See Appendix A for more of all these details.

The identity bijection $\id \in \pn$ is the natural alignment between $G$ and $H$.
When $\re=1$ we have that $G,H$ are almost surely isomorphic (via isomorphism $\id$), and when $\re=0$
we have that $G$ and $H$ are independent (i.e.~the indicators for all edges
of both graphs are collectively independent). If all Bernoulli parameters $p_{i,j}$
are equal to each other then $G$ and $H$ are Erdos-Renyi~random~graphs.

Associated with the Bernoulli parameters $\{ p_{i,j} \}_{\pairs \in \spairs}$,
denote their mean $$\mu:=\frac{1}{{n \choose 2}} \sum_{\pairs \in \spairs} p_{i,j}$$ and denote
their variance $$\sigma^2:=\frac{1}{{n \choose 2}} \sum_{\pairs \in \spairs} (p_{i,j}-\mu)^2.$$ We define
the {\it heterogeneity correlation}~$\rh$
\begin{equation}
\rh:=\frac{\sigma^2}{\mu (1-\mu)} .
\end{equation}
It is simple to show that $0 \leq \rh \leq 1$.
Furthermore, $\rh=0$ if and only if
all Bernoulli parameters $p_{i,j}$ are equal to each other (i.e.~$G$ and $H$ are
Erdos-Renyi random graphs), and $\rh=1$ if and only if all Bernoulli parameters are $0$ or $1$
(but, recall, the Bernoulli parameters are not all~$0$ and are not all $1$). See Appendix B for more details.
Note that $\rh$ is a measure of heterogeneity within $G$ (and within $H$) by virtue of its
numerator being the variance (a measure of spread) of the Bernoulli coefficients,
although this variance is normalized through division by the denominator of $\rh$,
where this denominator is a function of the global graph density.
(So, among distributions with a common global density, $\rh$ is just a multiple of~the~variance~$\sigma^2$.)

Note that edge correlation $\re$ is an inter-graph affect (between $G$ and $H$), whereas heterogeneity correlation
$\rh$ is an intra-graph affect. Unlike edge correlation $\re$, heterogeneity correlation $\rh$ is not a statistical
correlation. However, our results will demonstrate that $\rh$ is interchangeable with edge correlation $\re$ with
regard to creating alignment strength. We thus take the liberty of calling $\rh$ ``correlation," but we do so in a looser,
nonstatistical sense, with the meaning that it generates similarity between $G$ and $H$ just like edge correlation does.

Finally, define the {\it total correlation} $\rT$ such that $\rT$ satisfies
\begin{equation}
1-\rT:=(1-\re)(1-\rh).
\end{equation}

\section{Alignment strength is total correlation, asymptotically  \label{S:main}} 

In this section we state and prove our main theoretical result, Theorem \ref{thm:after}, that when $G,H$ have a correlated
Bernoulli distribution then the identity bijection $\id \in \pn$ (the natural
alignment here) has alignment strength asymptotically equal to the distribution's total correlation $\rT$.

Let $e_G$ and $e_H$ denote the number of edges in $G$ and $H$, respectively, and let $\dens_G:=\frac{e_G}{{n \choose 2}}$
and $\dens_H:=\frac{e_H}{{n \choose 2}}$ respectively denote the densities of $G$ and $H$.
\begin{lemma} \label{thm:first}
For any graphs $G$, $H$ on common vertex set $[n]$, and any $\phi \in \pn$, it holds that
$$\st (G,H,\phi)=1-\frac{ \frac{d(G,H,\phi)}{{n \choose 2}}  }{\dens_G \left ( 1-\dens_H \right )+  \left ( 1-\dens_G \right ) \dens_H}.$$
\end{lemma}

\noindent {\bf Proof:} With $G$ and $H$ fixed,
consider random $\varphi \in \pn$ with a discrete-uniform distribution; the expected value of $d(G,H,\varphi)$
is $ \frac{1}{n!}\sum_{\phi' \in \pn}d(G,H,\phi')$.
We next equivalently compute the expected value of $d(G,H,\varphi)$ using
linearity of expectation over the sum of its indicators in Equation~\ref{eqn:astr}. Observe that,
for any two vertices that form an edge in $G$, the probability that $\varphi$ maps them to a nonedge of
$H$ is $\frac{{n \choose 2}-e_H}{{n \choose 2}}$, and, for any two nonadjacent vertices of $G$, the probability
that $\varphi$ maps them to an edge of $H$ is $\frac{e_H}{{n \choose 2}}$; the expected value of $d(G,H,\varphi)$ is thus
\begin{eqnarray}
& &  \frac{1}{n!}\sum_{\phi' \in \pn}d(G,H,\phi') \nonumber \\
  &=& e_G \cdot \frac{{n \choose 2}-e_H}{{n \choose 2}}+ \left  ({n \choose 2}-e_G \right )
\cdot \frac{e_H}{{n \choose 2}} \nonumber \\  &=& {n \choose 2} \cdot \Big [  \dens_G \left ( 1-\dens_H \right )
+  \left ( 1-\dens_G \right ) \dens_H \Big ]. \label{eqn:b}
\end{eqnarray}
The desired result then follows from substituting Equation \ref{eqn:b} into the definition of $\st (G,H,\phi)$ in Equation~\ref{eqn:st}.~$\qed$\\

In the rest of this section we will state and prove limit results for random correlated Bernoulli graphs $G$, $H$.
This context requires us to consider a sequence of experiments ---for each value of $n=1,2,3,\ldots$ ---wherein
the chosen edge correlation $\re$ is a function of $n$, and the chosen Bernoulli
parameters $\{ p_{i,j} \}_{\pairs \in \spairs}$ are also functions of $n$, and thus $\rh$ and $\rT$
are also functions of $n$. For ease of notation, we do not explicitly write argument $n$ in these
functions. However, we will require that there exists a positive lower bound for $\mu$
over all $n$, and as well that there exists an upper bound less than $1$ for $\mu$ over all $n$. (Note that
since $\mu$ is a function of $n$, we have that the $\mu$ are a sequence, so
the following limit result is expressed as a difference that converges as stated, rather than convergence to $\mu$,
which would not make technical sense. Similarly for the other results here.)

\begin{lemma} \label{thm:lead1}
We have $ \dens_G -\mu \stackrel{a.s.}{\rightarrow} 0$ and $ \dens_H - \mu \stackrel{a.s.}{\rightarrow} 0$.
\end{lemma}

\noindent {\bf Proof:} Clearly $\e (\dens_G) =\mu$. Also, $e_G$ is the sum of ${n \choose 2}$ independent
Bernoulli random variables, and thus its variance is bounded by ${n \choose 2}$, thus the variance of
$\dens_G:=\frac{e_G}{{n \choose 2}}$ is of order O$(n^{-2})$. Next, by Chebyshev's Inequality, for any $\epsilon>0$, \
$\p \left [ \left | \dens_G- \mu  \right | \geq \epsilon   \right ]
\leq \frac{1}{\epsilon ^2} \textup{Var} \left (   \dens_G \right )$;
since this probability is O$(n^{-2})$ when $\epsilon$ is fixed, it has finite sum over $n=1,2,3,\ldots$.
Thus, since $\epsilon$ is arbitrary, by the Borel-Cantelli Theorem  $ \dens_G -\mu \stackrel{a.s.}{\rightarrow} 0$, as desired. $\qed$

\begin{theorem} \label{thm:lead2}
We have $\frac{d(G,H,\id)}{{n \choose 2}} - 2 (1-\re) \Big (  \mu(1-\mu) - \sigma^2 \Big ) \stackrel{a.s.}{\rightarrow} 0 $
\end{theorem}
\noindent {\bf Proof:}
We begin by taking the expected value of $d(G,H,\id)$;
\begin{eqnarray}
& & \e \Big [ d(G,H,\id) \Big ] \nonumber \\ & = & \e \left [   \sum_{\pairs \in \spairs} \q \Big ( \ \q \big (i \sim_G j \big )
\ \ne \ \q \big (i \sim_H j \big ) \ \Big )  \right ]  \nonumber \\
& = & \sum_{\pairs \in \spairs} 2 (1-\re) p_{i,j}(1-p_{i,j})\nonumber
\\  &=& 2 (1-\re) {n \choose 2}  \Big (  \mu(1-\mu) - \sigma^2 \Big ), \label{eqn:num}
\end{eqnarray}
thus $\e \Big [ \frac{d(G,H,\id)}{{n \choose 2}} \Big ]= 2 (1-\re) \Big (  \mu(1-\mu) - \sigma^2 \Big )$.\\
Next, $d(G,H,\id)$ is is the sum of ${n \choose 2}$ independent
Bernoulli random variables, and thus its variance is bounded by ${n \choose 2}$, thus the variance of
$\frac{d(G,H,\id)}{{n \choose 2}}$ is of order O$(n^{-2})$. Next, by~Chebyshev's~Inequality,~for~any~$\epsilon>0$, \
$\p \left [ \left | \frac{d(G,H,\id)}{{n \choose 2}} - 2 (1-\re) \Big (  \mu(1-\mu) - \sigma^2 \Big ) \right | \geq \epsilon   \right ]
\leq \frac{1}{\epsilon ^2} \textup{Var} \left (  \frac{d(G,H,\id)}{{n \choose 2}} \right )$;
since this probability is O$(n^{-2})$ when $\epsilon$ is fixed, it has finite sum over $n=1,2,3,\ldots$.
Thus, since $\epsilon$ is arbitrary, by the Borel-Cantelli Theorem $\frac{d(G,H,\id)}{{n \choose 2}} - 2 (1-\re)
 \Big (  \mu(1-\mu) - \sigma^2 \Big ) \stackrel{a.s.}{\rightarrow} 0 $, as desired. $\qed$\\

\noindent The following is the main result of this section, and is our main theoretical result.

\begin{theorem} \label{thm:after} It
holds that $\st (G,H, \id ) -\rT \stackrel{a.s.}{\rightarrow} 0 $
\end{theorem}
\noindent {\bf Proof:} By Lemma \ref{thm:lead1}, \
$ \dens_G -\mu \stackrel{a.s.}{\rightarrow} 0$ and $ \dens_H - \mu \stackrel{a.s.}{\rightarrow} 0$.
Because $\dens_G$, $\dens_H$ and $\mu$ are bounded, we thus have
that $\dens_G \left ( 1-\dens_H \right )+  \left ( 1-\dens_G \right ) \dens_H - 2 \mu(1-\mu) \stackrel{a.s.}{\rightarrow} 0$.
Now, by Theorem \ref{thm:lead2}, we have that
$\frac{d(G,H,\id)}{{n \choose 2}} - 2 (1-\re) \Big (  \mu(1-\mu) - \sigma^2 \Big ) \stackrel{a.s.}{\rightarrow} 0$; since the
relevant sequences are bounded, and $\mu$ is bounded away from $0$ and $1$, we have that
$$\frac{ \frac{d(G,H,\id)}{{n \choose 2}} }{\dens_G \left ( 1-\dens_H \right )+  \left ( 1-\dens_G \right ) \dens_H } \  - \
\frac{ 2 (1-\re) \Big (  \mu(1-\mu) - \sigma^2 \Big )}{2 \mu(1-\mu)}\stackrel{a.s.}{\rightarrow} 0.  $$
Applying Lemma \ref{thm:first} and the definitions of $\rh$ and $\rT$
we thus have from above that $(1-\st (G,H, \id ) )- (1 - \rT) \stackrel{a.s.}{\rightarrow} 0 $, which proves Theorem \ref{thm:after}. $\qed$

\section{Graph matchability and total correlation $\rT$  \label{S:matchability}}

In this section we empirically demonstrate  in broad families of parameter settings where $\re$\ and $\rh$ vary,
that success of an approximate seeded graph matching algorithm is a function~of~$\rT$.

Our setting is where $G$, $H$ are correlated Bernoulli graphs on vertex set $[n]$. The graph matching problem
is to compute $\phi_{GM} \in \arg \min_{\phi \in \pn}d(G,H,\phi)$. In the {\it seeded graph matching problem}, there
are $s$ {\it seeds}, without loss of generality they are the vertices $1,2,\ldots ,s$, and there are
$m:=n-s$ ambiguous vertices, which are the other vertices $s+1,s+2, \ldots, n$. The meaning of {\it seeded graph matching}
is that the feasible region $\phi \in \pn$ of the graph matching problem is restricted to $\phi \in \pn$ that satisfy $\phi (i)=i$ for all seeds
$i=1,2,\ldots,s$. The graphs $G$ and $H$ are separately observed and the identities of the ambiguous
vertices are unobserved for the optimization,
so that the natural alignment, which is the identity bijection $\id$, is only seen for the seeds. If
the seeded graph matching solution is $\id$ then we say that $G$ and $H$ are {\it matchable}.

Even a modest number of seeds can make a very significant increase in the likelihood
that $G$ and $H$ are matchable [\cite{sgm2}]. Our illustration in this section will be for realistically
sized graphs, on the order of a thousand vertices, and we utilize seeds because they
will be quite helpful in obtaining reasonable probability of matchability.
Unfortunately, exact graph matching --even seeded graph matching-- is intractable,
only solvable on the smallest, toy examples. So we utilize an approximate seeded graph matching algorithm;
the specific one we use is the SGM Algorithm [\cite{sgm1}, \cite{sgm2}], which has been demonstrated to have many
nice theoretical properties, and it is efficient and quite effective (see \cite{sgm1}, \cite{sgm2}, \cite{relax}). In this section,
we will say that $G$ and $H$ are matchable if the SGM-generated approximate seeded graph matching solution
is the identity bijection $\id$.

In the experiments that we will perform, we will sample $G$, $H$ from a correlated Bernoulli distribution
for different values of $\re$ and $\rh$; the values of the Bernoulli coefficients
$\{ p_{i,j} \}_{\pairs \in \spairs}$ are selected as follows, in order to obtain specified values of $\rh$.
Given any real number $p \in (0,1)$ and real number $\delta \in [0, \min \{ p, 1-p \} ]$, we independently randomly
sample $\{ p_{i,j} \}_{\pairs \in \spairs}$ from the uniform distribution on the interval $(p-\delta, p+\delta)$.
Note that the afore-defined Bernoulli parameter variance $\sigma^2$ has expected value~$\frac{\delta^2}{3}$, and $\sigma^2$ will be approximately
$\frac{\delta^2}{3}$ for large values of $n$. For a fixed $p$, as $\delta$ goes from $0$ to $\min \{ p, 1-p \}$,
the value of  $\rh=\frac{\sigma^2}{\mu (1-\mu)} \approx \frac{\delta^2}{3p(1-p)} $ monotonically increases
from $0$ to $\frac{1}{3} \cdot \frac{1-p}{p}$ if $p \geq \frac{1}{2}$ and $\frac{1}{3} \cdot \frac{p}{1-p}$ if $p \leq \frac{1}{2}$.
In this section, when we report values of $\re$ and $\rh$, we mean that we selected $\delta$ so that the
approximate value of $\rh$ is as reported.

\begin{figure}[h!]
	\centering
	\includegraphics[width=3.7in]{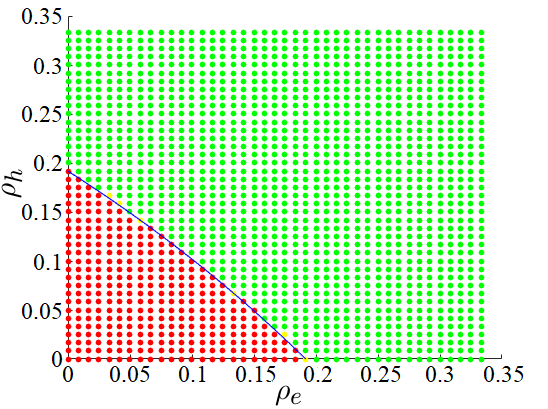}
    \caption{Matchability experiment for $m=850$, $s=150$, $p=\frac{1}{2}$.}		
	\label{fig:exper1}
\end{figure}

We did three batches of experiments. In the first batch of experiments, for each value of $\re = 0, \frac{1}{120},
\frac{2}{120}, \frac{3}{120}, \ldots, \frac{1}{3}$ and $\rh = 0, \frac{1}{120},
\frac{2}{120}, \frac{3}{120}, \ldots, \frac{1}{3}$, we did $60$ replicates of obtaining random graphs $G$, $H$ with
$m=850$ ambiguous vertices and $s=150$ seeds from a correlated Bernoulli distribution with edge correlation $\re$ and
heterogeneity correlation $\rh$ based on $p=\frac{1}{2}$, and we performed seeded graph matching with the SGM algorithm. If all $60$ replicates were
matchable then we plotted a green dot in Figure \ref{fig:exper1} at the appropriate coordinates, if between $1$ and $5$ of the $60$
replicates were not matchable then we plotted a yellow dot in the figure, and if more than $5$ of the $60$ replicates
were not matchable then we plotted a red dot. The blue curve in the figure is the set of all pairs of $\re$, $\rh$ such that
$\rT=\frac{23}{120}$. \

In these experiments and those below, the transition from matchable to anonymized (i.e., not matchable) occurs at a level set of $\rT$.
We note here that numerous results in the literature have studied this matchability phase transition as a function of edge correlation $\re$ (see, for example, [\cite{match1,match2,sgm2}]) and a few papers have considered the impact of network heterogeneity on matchability (see, for example, [\cite{match3,match4}]).
In the parameterized correlated Bernoulli distribution considered above, these empirical results novelly suggest the form by which matchability is impacted by within and across graph correlation structure.
Further understanding this phase transition as a function of $\rT$ is a necessary next step to understand the dual roles that graph structure ($\rh$) and graph pairedness ($\re$) play in network alignment problems both theoretical and practical.

\begin{figure}[h!]
	\centering
	\includegraphics[width=3.7in]{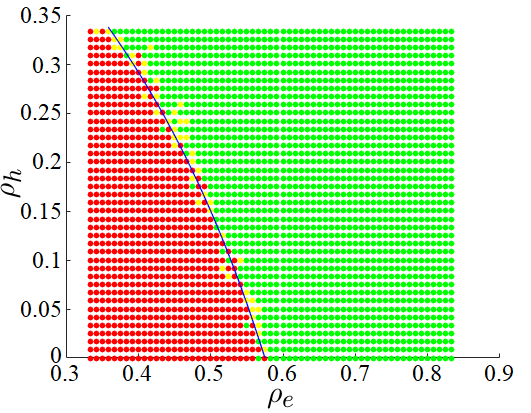}
    \caption{Matchability experiment for $m=850$, $s=9$, $p=\frac{1}{2}$}		
	\label{fig:exper2}
\end{figure}

The second batch of experiments differed just in that there were only $s=9$ seeds (with $m=850$ as before), and the range of values of $\re$ was
$\frac{1}{3}$ to  $\frac{5}{6}$ in increments of $\frac{1}{120}$; the results are similarly displayed in Figure \ref{fig:exper2}, and the blue curve
in the figure is the set of all pairs of $\re$, $\rh$ such that
$\rT=\frac{69}{120}$. In these experiments, we again see the transition in matchability
at a level set of $\rT$, although the transition is looser due to fewer seeds being considered in this problem setup.

\begin{figure}[h!]
	\centering
	\includegraphics[width=3.7in]{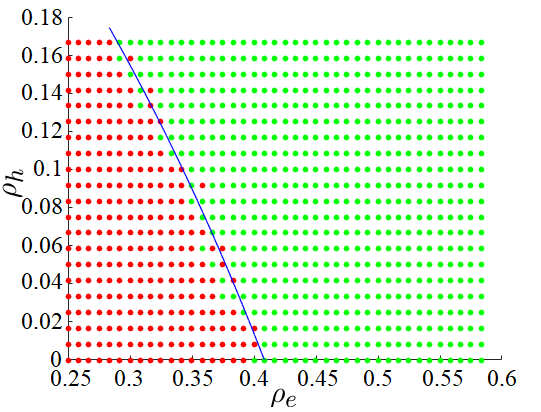}
    \caption{Matchability experiment for $m=850$, $s=22$, $p=\frac{1}{3}$}		
	\label{fig:exper3}
\end{figure}

The third batch of experiments differed just in that there were $s=22$ seeds, and now $p=\frac{1}{3}$, the range of values of $\re$
was $\frac{1}{4}$ to $\frac{7}{12}$ in increments of $\frac{1}{120}$, and the range of values of $\rh$ was
$0$ to $\frac{1}{6}$ in increments of $\frac{1}{120}$; the results are similarly displayed in Figure \ref{fig:exper3},
and the blue curve in the figure is the set of all pairs of $\re$, $\rh$ such that
$\rT=\frac{49}{120}$.
In these experiments, we again see the transition in matchability
at a level set of $\rT$.

We then repeated the above experiments for each combination
of: total number of vertices $300$ or $600$, number of seeds
seeds $5\%$ or $10\%$ of the vertices, and values of $p$ being
$\frac{1}{2}$ or $\frac{1}{3}$. In all eight such combinations
the result of the experiments were like the above; namely, matchability was a function of $\rT$.

Note that matchability is not universally a function of just $\rT$.
For example, the number of vertices and the number of seeds have a dramatic affect on matchability.
The empirical
demonstrations in this section of matchability as a function of $\rT$ are limited to
families of correlated Bernoulli distribution parameterizations of the type that we have used here. New work
will be needed to obtain theorems that universally and fully account for matchability.
But, nonetheless,
we have empirically demonstrated in broad families of parameter settings that the phase transition in matchability
occurs at a level set of $\rT$, which supports the importance and utility of $\rT$ as a meaningful measure of graph correlation.

\section{Graph matching runtime and total correlation $\rT$  \label{S:runtime}}

Similar to the previous section, in this section we empirically demonstrate,  in broad families of parameter settings where
$\re$\ and $\rh$ vary, that the running time of exact seeded graph matching via binary integer linear
programming is~a~function~of~$\rT$.

We consider exact seeded graph matching here because the approximate seeded graph matching algorithms have running times
that are relatively monolithic (when the number of vertices are fixed)
and not sensitive enough to the parameters in the random graph distribution.
Unfortunately, exact graph matching is intractable [\cite{thirtyyears}], and can only be done
for small examples; we will work with graphs that have $20$ ambiguous vertices.

For this section, the random graphs $G$,$H$ have correlated Bernoulli distributions, for various
values of $\re$ and $\rh$. The Bernoulli parameters are chosen in exactly the manner of the
previous section, Section \ref{S:matchability}; there is a fixed value $p$, and then $\delta$ are
selected to attain desired values of $\rh$ in the manner described in the previous section.

We next formulate the binary integer linear program for seeded graph matching.
For graphs $G$ and $H$, say their adjacency matrices are $A$ and $B$,
respectively, and say that there are $s$ seeds and $m$ ambiguous vertices. We partition
$A = \bigl [\begin{smallmatrix} A_{11} & A_{12} \\ A_{21} & A_{22} \end{smallmatrix} \bigr ]$ and
$B = \bigl [\begin{smallmatrix} B_{11} & B_{12} \\ B_{21} & B_{22} \end{smallmatrix} \bigr ]$, where
$A_{11},B_{11}~\in~\{ 0,1 \}^{s \times s}$,  $A_{12},B_{12}~\in~\{ 0,1 \}^{s \times m}$,
 $A_{21},B_{21}~\in~\{ 0,1 \}^{m \times s}$, and $A_{22},B_{22} \in \{ 0,1 \}^{m \times m}$.
(Note that $A_{12}=A_{21}^T$ and $B_{12}=B_{21}^T$ here, since $A$ and $B$ are symmetric, but we do not use this
fact in the formulation below so that the formulation is expressed even more generally.)
Let $I$ denote the
 identity matrix (subscripted with its number of rows and columns), let $0$ subscripted denote the matrix
 of zeros of subscripted size, let $\vec{1}$ denote the column vector of ones with subscripted number of entries,
 let $\vec{0}$ denote the column vector of zeros with subscripted number of entries,
 let $\otimes$ denote the Kronecker product of matrices,
 let $\| \cdot \|_1$ denote the $\ell_1$ vector norm for matrices (this norm is evaluated
 by taking the sum of absolute values of the matrix entries), for any matrix $N$ let $\vc N$ denote the column vector
 which is the concatenation of the columns of $N$ (first column of~$N$, then second column of $N$, etc., then last column of $N$),
 and let $\ppm$ denote the set of
 $m \times m$ permutation matrices. Clearly, the seeded graph matching problem is $\min_{P \in \ppm}
 \| A - \bigl [\begin{smallmatrix} I_{s \times s} & 0_{s \times m} \\ 0_{m \times s} & P \end{smallmatrix} \bigr ]
 B \bigl [\begin{smallmatrix} I_{s \times s} & 0_{s \times m} \\ 0_{m \times s} & P \end{smallmatrix} \bigr ]^T \|_1$.
 By permuting columns of the matrix in the norm, we get an equivalent formulation of the seeded graph matching
 problem as:
 \begin{align*}
 \min_{P \in \ppm}
 \| A \bigl [\begin{smallmatrix} I_{s \times s} & 0_{s \times m} \\ 0_{m \times s} & P \end{smallmatrix} \bigr ]
 - \bigl [\begin{smallmatrix} I_{s \times s} & 0_{s \times m} \\ 0_{m \times s} & P \end{smallmatrix} \bigr ] B  \|_1.
 \end{align*}
 Expanding this, we get an equivalent formulation of the seeded graph matching problem as:
 \begin{equation} \label{eqn:sgmprob} \hspace{-.4in}
 \min_{P \in \ppm}
 \Big ( \| A_{12}P-B_{12} \|_1 + \| A_{21}-PB_{21}\|_1 + \| A_{22}P-PB_{22} \|_1  \Big ) .
 \end{equation}
Now, because of the absolute values in $\| \cdot \|_1$, we add artificial  variables to obtain simple linearity. For example,
(just) minimizing  $\| A_{22}P-PB_{22} \|_1$ subject to $P \in \ppm$ is equivalent to minimizing the
sum of the entries of $E+E'$ subject to $A_{22}P-PB_{22}+E-E'=0_{m \times m}$, \ \ $P \in \ppm$,  \ \  $E,E' \in \{ 0,1 \} ^{m \times m}$.
Of course, there are additional $\| \cdot \|_1$ terms in the objective function in Equation~\ref{eqn:sgmprob},
but the same approach can be used, so that seeded graph matching is equivalent to
\begin{align*}
\min  \, \, &\bigl [\begin{smallmatrix} \vec{0}_{m^2}  \\
\vec{1}_{2m^2+4ms} \end{smallmatrix} \bigr ]^Tx\\
\text{ s.t. }&[M|E] x=b\\
&x \in \{ 0,1 \}^{3m^2+4ms}
\end{align*}
where the first $m^2$ entries of $x$ are $\vc P$, and $M$ and $E$ and $b$ are given by:
\[
M=
\left [ \begin{array}{c}
I_{m \times m} \otimes A_{22} - B_{22}^T \otimes I_{m \times m} \\
I_{m \times m} \otimes A_{12}  \\  B_{21}^T \otimes I_{m \times m} \\
I_{m \times m} \otimes \vec{1}_m^T \\  \vec{1}_m^T \otimes I_{m \times m}
\end{array} \right ]
\]
\[
E=
\left [ \begin{array}{cc}
I_{(m^2+2ms) \times (m^2+2ms)}  & -I_{(m^2+2ms) \times (m^2+2ms)} \\
0_{      2m \times  (m^2+2ms)}  & 0_{2m \times (m^2+2ms)}
\end{array} \right ]
\]
\[
b=\left [ \begin{array}{c} \vec{0}_{m^2} \\ \vc B_{12} \\ \vc A_{21} \\ \vec{1}_m \\ \vec{1}_m \end{array} \right ]
\]
We solve the above binary integer linear program exactly using the optimization package GUROBI.
The yardstick for runtime that we have chosen to adopt
is the number of simplex iterations performed by GUROBI;
this measure has the advantage of reducing many sources of platform variability.

%
%
%
%

\begin{figure}[h!]
	\centering
	\includegraphics[width=5.6in]{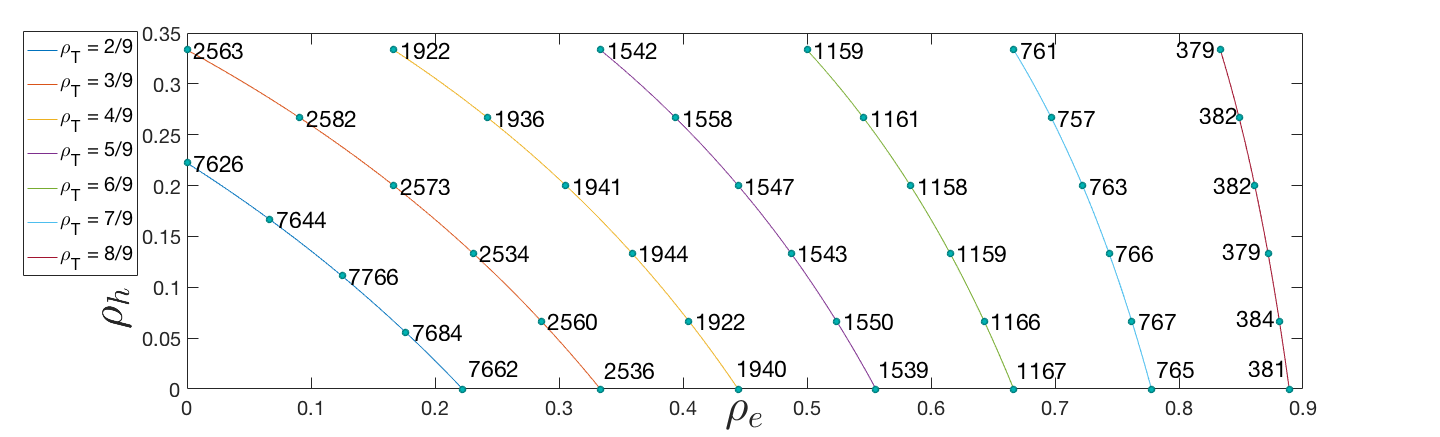}
    \caption{Runtime experiment for $m=20$, $s=480$, $p=\frac{1}{2}$.}		
	\label{fig:exper4}
\end{figure}
\begin{figure}[h!]
	\centering
	\includegraphics[width=5.6in]{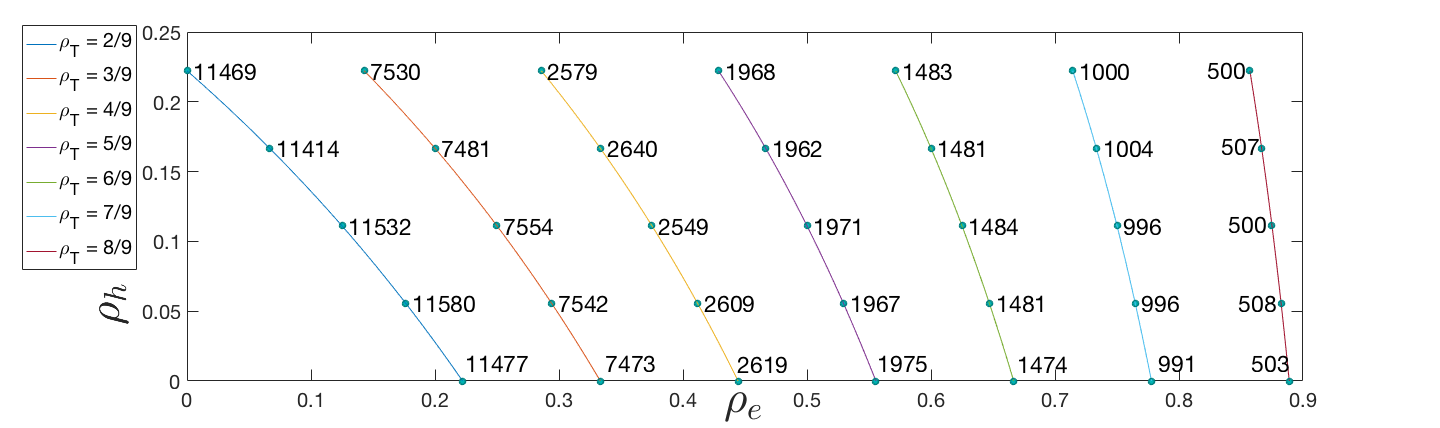}
    \caption{Runtime experiment for $m=20$, $s=480$, $p=\frac{3}{5}$.}		
	\label{fig:exper5}
\end{figure}
\begin{figure}[h!]
	\centering
	\includegraphics[width=5.6in]{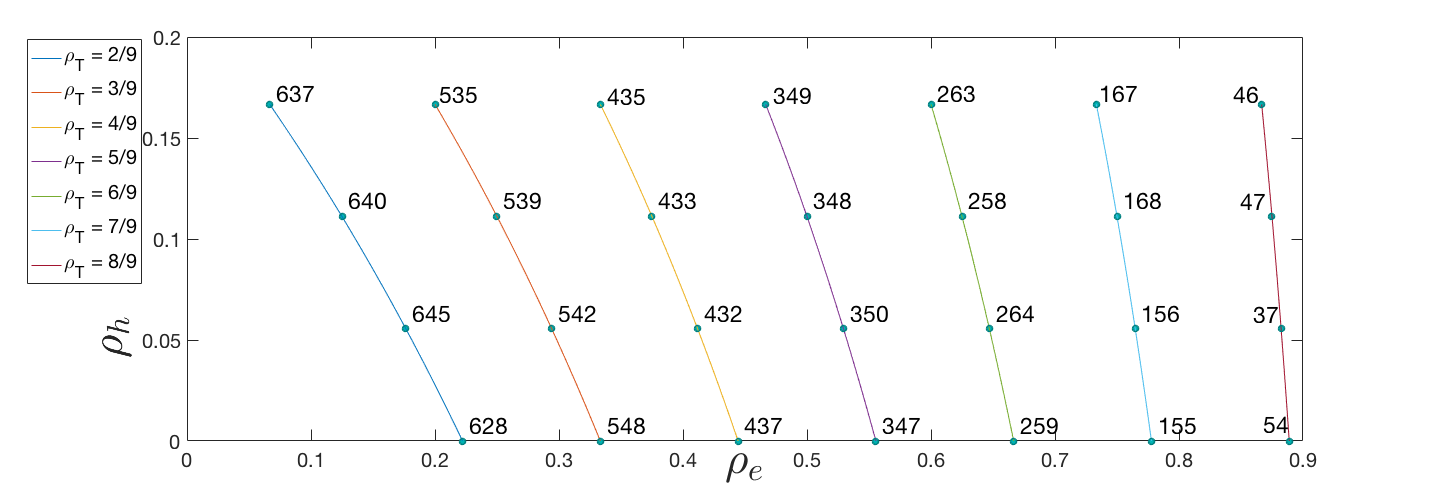}
    \caption{Runtime experiment for $m=20$, $s=480$, $p=\frac{1}{3}$.}		
	\label{fig:exper6}
\end{figure}

We performed three batches of experiments. In the first batch of experiments, for each value~of $\rT = \frac{2}{9},
\frac{3}{9}, \frac{4}{9}, \ldots, \frac{8}{9}$, we selected various pairs of
$\re$, $\rh$ which have $1-\rT=(1-\re)(1-\rh)$ for the given value of $\rT$; the values of $\rh$
are achieved based on $p=\frac{1}{2}$, and the chosen pairs $\re$, $\rh$ are the points plotted with a dot in Figure~\ref{fig:exper4}.
For each such pair $\re$, $\rh$ we  did $60$ replicates of obtaining random graphs $G$, $H$ with
$m=20$ ambiguous vertices and $s=480$ seeds from a correlated Bernoulli distribution with edge correlation $\re$ and
heterogeneity correlation $\rh$, and we solved the seeded graph matching problem for $G$, $H$ exactly using GUROBI.
The average runtimes (measured by the number of simplex iterations performed by GUROBI)
are printed above each pair $\re$, $\rh$ at the appropriate coordinates in Figure \ref{fig:exper4}.
The smooth curves on the plot are~the~level~sets~of~$\rT$.

These experiments, and those below, suggest that in this parametrized Bernoulli graph model the algorithmic runtimes are approximately constant on the level sets of $\rT$.
The results in Section~\ref{S:matchability} suggest that the phase transition of matchability occurs at a level set of $\rT$, and these results further reinforce the novel overarching notion:
that the theoretic and algorithmic difficulty of matching is a function of $\re$ and $\rh$ only through $\rT$.
Alone, $\re$ and $\rh$ are insufficient to capture this theoretic and algorithmic difficulty.

The second and third batch of experiments are exactly like the first batch, except that
for the second batch of experiments the values of $\rh$ are based on $p=\frac{3}{5}$ and
the results are displayed in Figure \ref{fig:exper5}, and
for the third batch of experiments the values of $\rh$ are based on $p=\frac{1}{3}$ and
the results are displayed in Figure \ref{fig:exper6}.
Note that the ranges of $\rh$ are different in Figures \ref{fig:exper4}, \ref{fig:exper5}, and \ref{fig:exper6}
because different values of $p$ put different limitations~on~$\delta$.

Just like for matchability in the previous section, it must be pointed
out that the runtime of exact seeded graph matching via binary integer
linear programming is not universally a function of $\rT$.
Of course, the number of vertices ---particularly the number of ambiguous vertices---
has a dominant role in the runtime, and the above experiments show that the graph density
likewise plays a very large role. Nonetheless, for families of correlated Bernoulli distributed
graphs similar to the ones in the experiments above, we see  within a family that the runtime
is a function of $\rT$.

\section{Discussion and future work \label{S:discussion}}

The correlated Bernoulli random graph model considered herein contains many important families
of random graph models as subfamilies including stochastic blockmodels [\cite{sbm,airoldi2008mixed}],
random dot product graphs [\cite{young2007random,athreya2017statistical}], and more general latent
position random graph [\cite{hoff2002latent}]. While the edge independent assumption inherent to these
models is often not satisfied in real data applications, nonetheless (conditionally) edge-independent
random models have shown great utility in capturing statistically relevant structure in a host of real data
applications from modeling connectomic
structure [\cite{priebe2017semiparametric,lyzinski2017community,levin2017central}], to capturing
community and user-level behavior in social networks [\cite{xu2014dynamic,patsolic2017vertex}].
Moreover, these models provide a theoretically tractable environment in which to explore important
statistical concepts such as estimation
consistency [\cite{bickel2009nonparametric,rohe2011spectral,sussman2014consistent}],
consistent hypothesis testing [\cite{tang2014nonparametric,MT2,lei2014goodness}], and network
de-anonymization [\cite{eralg2,match1}], among others. Indeed, it is this appealing
mix of theoretical tractability and practical utility that have made these graph models an increasingly popular
option in the statistical network inference community.

In this paper we prove in a very broad random graph setting
---specifically, when
$G$ and $H$ have a correlated Bernoulli distribution---
that the alignment
strength of the natural $G$, $H$ alignment is asymptotically equal to the total
correlation $\rT$ in the distribution. After this, we empirically demonstrate, for
types of families within the distribution, that both matchability and
exact-solution-runtime for seeded graph matching of $G$,~$H$ are functions of the total correlation $\rT$.

Graph matching and seeded graph matching are extremely important in many disciplines; see the surveys [\cite{thirtyyears}] and
[\cite{tenyears}].
Unfortunately, these problems are intractable; in their full generality they are NP-hard.
Obtaining a function of the
distribution parameters that universally predicts matchability via approximate algorithms would be a huge advance in
theoretical understanding and in practice. Likewise, it would be a huge advance to predict
exact-solution-runtime from a function of the distribution parameters, and it would not just be the number of
vertices---the other parameters play a large role. The goals of obtaining these universal
functions has not been achieved here; the families we use here are general but not universal.
But a universal result will include our families as special cases, thus $\rT$ will play an
important role.

There are a number of matchability results already known, see [\cite{match1,match2,match3,sgm2,relax,match4}]. However, for the most part these are
asymptotic results that do not specify the particular constants involved,
and leave gaps in the parameter possibilities where the results are
silent.
In particular, the empirical matchability demonstrations in this paper are not
predictable from the previously known matchability asymptotics.
 Many of the known matchability results explicitly or implicitly involve edge correlation $\re$. The formulation of
$\rh$ is new to this paper, and $\rT$ is also new to this paper. Thus we are now opening
a fertile new avenue for proof-of-matchability results based on $\rh$ and $\rT$, in the spirit
of the existing results for $\re$ and also for more powerful types of results.

\section*{Acknowledgments}

The authors are grateful to
The Maryland Advanced Research Computing Center for use of their supercomputer
to conduct the computational experiments. An anonymous contributor
made a very useful observation which led to streamlining the main result's proof.
The referees' and editor's feedback and remarks  greatly
strengthened this article, and are very much appreciated.
Our research was sponsored by the Air Force Research Laboratory and DARPA,
under agreement numbers FA8750-18-2-0035 and FA8750-17-2-0112.
The U.S. Government is authorized to reproduce and distribute reprints for Governmental purposes notwithstanding any copyright notation thereon.The views and conclusions contained herein are those of the authors and should not be interpreted as representing official policies or endorsements, expressed or implied, of Air Force Research Laboratory, DARPA, or the U.S. Government.

\bibliographystyle{plainnat}
\bibliography{refs}

\section*{Appendix}

We here provide some details about correlated Bernoulli random graphs. Notation here is as defined in the article.\\

\noindent {\bf Section A:} For any $\pairs \in \spairs$ such that $0<p_{i,j}<1$, suppose that
$\q (i \sim_G j)$ is a Bernoulli$(p_{i,j})$ random variable and
$\q (i \sim_H j)$ is a  Bernoulli$(p_{i,j})$ random variable, and suppose that
the two random variables $\q (i \sim_G j)$ and $\q (i \sim_H j)$ have Pearson correlation coefficient $\re$;
we derive the joint distribution of $\q (i \sim_G j)$ and $\q (i \sim_H j)$ as follows:
\begin{eqnarray*}
\re & = & \frac{\textup{Cov}\left [ \q (i \sim_G j),\q (i \sim_H j) \right ]}{\sqrt{\textup{Var}[\q (i \sim_G j)]},\sqrt{\textup{Var}[\q (i \sim_H j)]}}\\
& = & \frac{ \e [\q (i \sim_G j)\q (i \sim_H j)]- \e [\q (i \sim_G j )] \cdot \e [\q (i \sim_H j)] }{ \sqrt{p_{i,j}(1-p_{i,j})} \sqrt{p_{i,j}(1-p_{i,j})} }\\
& = & \frac{\p [  i \sim_G j \ \& \ i \sim_H j ] -p_{i,j}^2}{p_{i,j}(1-p_{i,j})},
\end{eqnarray*}
from which we obtain $\p [  i \sim_G j \ \& \ i \sim_H j ]= p_{i,j}^2 + \re p_{i,j}(1-p_{i,j}) $. Because
$\q (i \sim_G j)$ and $\q (i \sim_H j)$ are each marginally Bernoulli$(p_{i,j})$, we obtain
that $\p [  i \sim_G j \ \& \ i \not \sim_H j ]=\p [  i \not \sim_G j \ \& \ i \sim_H j ]=p_{i,j}-
\Big ( p_{i,j}^2 + \re p_{i,j}(1-p_{i,j}) \Big ) = (1-\re)p_{i,j}(1-p_{i,j})$, and also that
 $\p [  i \not \sim_G j \ \& \ i \not \sim_H j ]= (1-p_{i,j})- (1-\re)p_{i,j}(1-p_{i,j})=
 (1-p_{i,j})^2  + \re p_{i,j}(1-p_{i,j})$.

  Importantly, note that the joint distribution
 of  $\q (i \sim_G j)$ and $\q (i \sim_H j)$  is uniquely determined by $\re$. Also note that
 $\p [  \q (i \sim_G j) \ne \q (i \sim_H j)  ]=2(1-\re)p_{i,j}(1-p_{i,j})$. Also note that,
 conditioned on $\q (i \sim_G j)$, the random variable  Bernoulli$(\re \cdot \q (i \sim_G j) +
(1-\re) \cdot p_{i,j})$ results in the joint distribution above, which justifies the method in the article
of sampling $\q (i \sim_G j)$ and $\q (i \sim_H j)$ with marginal Bernoulli$(p_{i,j})$ distribution and
Pearson correlation coefficient $\re$. $\qed$\\

\noindent {\bf Section B:} We show that $\rh \leq 1$, with equality holding if and only if,
for all $\pairs \in \spairs$, it holds that $p_{i,j}$ is $0$ or $1$. Indeed,
\begin{eqnarray*}
& & 1-\rh \\
& = & 1-\frac{\sigma^2}{\mu (1-\mu)}\\
& = & \frac{            \mu (1-\mu)  -   \left (     \frac{   \sum_{\pairs \in \spairs} p_{i,j}^2 }{{n \choose 2}}     -\mu^2  \right )    }{\mu (1-\mu)}\\
& = & \frac{\sum_{\pairs \in \spairs}(p_{i,j}-p_{i,j}^2) }{{n \choose 2}\mu (1-\mu)}
\end{eqnarray*}
is clearly nonnegative and equals $0$ if and only if, for all $\pairs \in \spairs$ it holds that
$p_{i,j}=p_{i.j}^2$, i.e. it holds that $p_{i,j}$ is $0$ or $1$.
Thus $\rh \leq 1$ with equality holding if and only if,
for all $\pairs \in \spairs$, it holds that $p_{i,j}$ is $0$ or $1$.
(Except, recall, the Bernoulli parameters are not all~$0$ and are not all $1$, since $\rh$ would then not be defined.)
$\qed$

\end{document}